\documentclass[11pt,bezier]{article}
\usepackage{amsmath,amssymb,amsfonts,euscript,graphicx}

\title{{\LARGE\bf  Modules Whose Classical Prime Submodules Are Intersections of Maximal Submodules\thanks {The research of
 the second  author was in part supported by a grant from IPM (No. 89160031)}}}
\author{{\normalsize  {\bf M. Arabi-Kakavand${^{\rm a}}$ and   M.  Behboodi${^{\rm b}}$}}\thanks
{ Corresponding author}\\
  {\small{ ${^{\rm a,b}}$Department of Mathematical Sciences, Isfahan University of Technology}}\vspace{-1mm}\\ {\small{
Isfahan,
Iran, 84156-83111}}\\
{\small{$^{\rm b}$School of Mathematics,  Institute for Research in Fundamental Sciences (IPM), }} \\
 {\small{Tehran  Iran, 19395-5746}}\\
\footnotesize{${^{\rm a}}$m.arabi@math.iut.ac.ir}\vspace{-1mm}\\
\footnotesize{${^{\rm b}}$mbehbood@cc.iut.ac.ir}}
\date{}
\begin{document}
\maketitle
\begin{abstract}
{\noindent Commutative rings in which every prime ideal is the
intersection of maximal ideals are called Hilbert (or Jacobson)
rings. We propose to define classical Hilbert modules by the
property that {\it classical  prime} submodules are the
intersection of maximal submodules. It is shown that  all
co-semisimple modules as well as all Artinian modules are
classical Hilbert modules.  Also, every module over  a
zero-dimensional ring is classical Hilbert. Results illustrating
connections amongst the notions of classical Hilbert module and
Hilbert ring are also provided. Rings $R$ over which all
$R$-modules are classical Hilbert  are characterized. Furthermore,
we determine the Noetherian rings  $R$ for which all finitely
generated $R$-modules are classical Hilbert.\vspace{3mm}
  }
\end{abstract}

\section{Introduction}
   All rings in this paper are associative commutative
with identity $1\neq 0$  and modules are unital.
   Let  $M$ be  an $R$-module. If $N$ is a submodule (resp. proper submodule) of $M$, we write
$N\leq M$ (resp. $N<M$). The ideal $\{r\in R : rM\subseteq N\}$
will be denoted by  $(N:M)$. We call $M$ faithful if $(0 : M)=0$.
Also, we denote the classical Krull dimension of $R$ by dim$(R)$
and the  Jacobson radical of $R$ by $J(R)$.

A commutative ring $R$ is called a {\it Hilbert ring}, also {\it
Jacobson} or {\it Jacobson Hilbert ring}, if every prime ideal of
$R$ is the intersection of maximal ideals. This is obviously
equivalent to requiring that in each factor ring of R, the
nilradical coincides with the Jacobson radical.  The main interest
in Hilbert rings in commutative algebra and algebraic geometry is
their relation with Hilbert's Nullstellensatz; that is, if $R$ is
a Hilbert ring, then the polynomial ring $R[x_1, . . . , x_n]$ is
also a Hilbert ring (see for example \cite{AM,FP,FI,G,RJ}). This
notion was extended to noncommutative rings in several different
ways; see \cite{KW,P,W1,W2}.

In the literature, there are many different generalizations of the
notion of prime  ideals to  modules. For instance,  a proper
submodule $P$ of $M$ is called a {\it  prime submodule}   if
$am\in P$ for  $a\in R$ and $m\in M$ implies that  $m\in P$ or
$a\in (P:M)$.
 Prime submodules of modules were introduced  by J. Dauns \cite{D}  and have been studied
intensively since then (see for example \cite{AZ1,B2,B4}). Also, a
proper submodule $P$ of $M$ is called a {\it classical prime
submodule} if $abm\in P$ for $a$, $b\in R$ and $m\in M$ implies
that $am\in P$ or $bm\in P$. This notion of classical prime
submodule
 has been extensively studied by  the first  author in \cite{B1,B3}; see also \cite{BB,BN,BS}. Furthermore,
 in \cite{AZ1,B5,BK},  the authors
   use the terminology ``{\it weakly prime}'' to mean ``{\it classical prime}''.

There is already a generalization of the notion of commutative
Hilbert rings to modules. In fact, the notion of Hilbert modules
was introduced by Maani Shirazi and Sharif \cite{MSH}, by
requiring  that {\it prime submodules}  are  intersections of
maximal submodules. In this article we extend the notion of
commutative Hilbert rings to modules via classical prime
submodules. An $R$-module $M$ is a {\it classical Hilbert module}
(or simply {\it cl.Hilbert module}) if every classical prime
submodule of $M$ is an intersection of maximal submodules. In
Section 2, we study some properties of cl.Hilbert modules. Any
cl.Hilbert module is a Hilbert module but the converse need not be
true (see Example 2.1). It is shown that an $R$-module $M$ is a
cl.Hilbert module if and only if every non-maximal classical prime
submodule of $M$ is an intersection of properly larger classical
prime submodules (Theorem 2.5).  Any homomorphic image of a
cl.Hilbert module is a cl.Hilbert module (Proposition 2.6). This
yields that if $\oplus_{i\in I}M_{i}$ is a cl.Hilbert module, then
each $M_{i}$ $(i\in I)$ is a cl.Hilbert module (Corollary 2.8),
but the converse need not be true (see Example 2.9). Let $R$ be a
domain and $M$ be a cl.Hilbert $R$-module. If $N$ is any submodule
of $M$ such that $M/N$ is a torsion-free $R$-module, then $N$ is
also a cl.Hilbert $R$-module (see Proposition 2.13). This yields
the if $M$ is a cl.Hilbert module  over a domain $R$, then the
torsion submodule $T(M)$  is always a cl.Hilbert module. Moreover,
if $M$ is also torsion-free, then any pure submodule of $M$ is
also a cl.Hilbert module (see Corollary 2.14).
 It shown that all Artinian modules as well as all co-semisimple modules
 are
cl.Hilbert modules (see Example 2.2 (2) and Proposition 2.17 (3)).
Any torsion module over a one-dimensional domain is a cl.Hilbert
module (see Proposition 2.17 (2)). Also, it is shown that all
$R$-modules are cl.Hilbert if and only if ${\rm{dim}}(R)=0$ (see
Theorem 2.18). In Section 3 we investigate rings $R$ over which
every  finitely generated $R$-module is a cl.Hilbert module. In
particular, in  Theorem 3.6, we show that if $R$ is
 a  Noetherian  domain, then the following statements
are equivalent:\vspace{1mm}\\
(1)  Every finitely generated $R$-module is a
cl.Hilbert module.\\
(2)  The free $R$-module $R\oplus R$ is a
cl.Hilbert module.\\
(3) $R$ is both  a Hilbert ring and a  Dedekind domain.\\
(4) $R$ is a  Dedekind domain with $J(R)=0$.\\
 (5)  $R$ is either a field or a  Dedekind domain with infinity many  maximal ideals.

Furthermore,  we also characterize Noetherian rings  $R$ for which
every finitely generated $R$-module is a cl.Hilbert module (see
Theorem 3.7).

\section{Some properties of cl.Hilbert modules}

\hspace{3mm} Let $M$ be an $R$-module. Clearly every prime
submodule of $M$ is a classical prime submodule and, in case
$M=R$, where $R$ is any commutative ring, classical prime
 submodules and prime submodules coincide with prime ideals. But we may have a submodule $N$ in a
module $M$ that  is a classical prime submodule of $M$ but  is not
a prime submodule. In fact, if $R$ is a domain and $P$ is a
nonzero prime ideal in $R$, it is trivial to see that $P\oplus
(0)$,  $(0)\oplus P$ and $P(1, 1)$ are classical prime submodules
in the free module $M= R\oplus R$, but these are not prime
submodules (see also \cite[Example 3]{BK}). Thus  any
 cl.Hilbert module is a Hilbert module but  the following example shows
that the converse need not be true.\vspace{-1mm}\\

\noindent{\bf Example 2.1.} Let $R={\Bbb{Z}}[x]$. Since $R$ is a
Hilbert ring, by \cite[Proposition 2.9]{MSH},  the free
${\Bbb{Z}}[x]$-module ${\Bbb{Z}}[x]\oplus{\Bbb{Z}}[x]$ is a
Hilbert module. Now, for a prime number $p$ we put
$P=p{\Bbb{Z}}[x]+x{\Bbb{Z}}[x]$, which is the maximal ideal of
${\Bbb{Z}}[x]$ generated by the elements $p$ and $x$. We claim
that $P(p ,x)$ is a classical prime submodule of the free
${\Bbb{Z}}[x]$-module ${\Bbb{Z}}[x]\oplus{\Bbb{Z}}[x]$.
   To see this, let
$rs(f,g)\in P(p ,x)$, where $(f,g)\in {\Bbb{Z}}[x]\oplus
{\Bbb{Z}}[x]\setminus P(p ,x)$ and $r, s \in {\Bbb{Z}}[x]$. There
exists $z\in P$ such that $rs(f,g)=z(p, x)$, which implies that
$rsf=zp$ and $rsg=zx$.
   Suppose that $rs\neq 0$. Then any prime element $q$ of ${\Bbb{Z}}[x]$ which divides $rs$ must divide $z$, because $p$ and $x$ are co-prime in ${\Bbb{Z}}[x]$.
   It follows that $rs$ divides $z$. Hence there exists $z_1\in R$ such that $z=rsz_1$. This implies
    $f=z_1p$ and $g=z_1x$, i.e., $(f,g)=z_1(p,x)$. It follows that $z_1\not\in P$ since  $(f,g)\not\in
    P(p,x)$. But we have $rsz_1=z\in P$ and so $rs\in P$. Thus,  we have $r\in P$ or $s\in P$, which means that either
      $r(f,g)=rz_1(p,x)\in P(p,x)$ or $s(f,g)=sz_1(p,x)\in
    P(p,x)$. Thus $P(p,x)$ is a classical prime submodule of the free ${\Bbb{Z}}[x]$-module ${\Bbb{Z}}[x]\oplus{\Bbb{Z}}[x]$.
    Now we claim that  $P(p,x)$
     is not an intersection of maximal  submodules of ${\Bbb{Z}}[x]\oplus{\Bbb{Z}}[x]$.
     To see this, let $N$ be a maximal submodule of ${\Bbb{Z}}[x]\oplus{\Bbb{Z}}[x]$ such that
     $P(p,x)\subseteq N$. Since $N$ is a prime submodule, either  $(p,x)\in N$ or $P({\Bbb{Z}}[x]\oplus{\Bbb{Z}}[x])\subseteq
     N$. Since  $p, x\in P$, it follows that $(p,x)=p(1,0)+x(0,1)\in P({\Bbb{Z}}[x]\oplus{\Bbb{Z}}[x])$, which means that in any case, $(p,x)\in N$. Now, if  $P(p,x)$ is an intersection of maximal  submodules, then  we must have  $(p,x)\in
     P(p,x)$. It  follows that  $1\in P$, which  is a
     contradiction. Thus ${\Bbb{Z}}[x]\oplus{\Bbb{Z}}[x]$ is a
Hilbert $R$-module but it is not a cl.Hilbert $R$-module.

  We recall that if $U$, $M$
are $R$-modules, then following Azumaya, $U$ is called
$M$-injective if for any submodule $N$ of $M$, each homomorphism
$N\longrightarrow U$ can be extended to $M \longrightarrow U$ and
an $R$-module $M$ is called {\it co-semisimple} if every simple
module is $M$-injective  (see for example \cite[Chap. 4, Sec.
23]{Wis}). Also an $R$-module $M$ is called  {\it semisimple} if
$M$ is the direct sum of  all simple submodules.  Every semisimple
module is of course co-semisimple (see \cite[Proposition
23.1]{Wis}).

Next, we give several examples of cl.Hilbert modules. In
particular,  Parts (2) and (3)  of the following example show that
co-semisimple modules as well as all Artinian modules are
classical Hilbert modules.

\noindent{\bf Example 2.2}\vspace{2mm}\\
(1) Every Hilbert ring $R$ is a cl.Hilbert $R$-module (since every
classical prime submodule \indent of  $R$
is a prime ideal of $R$).\vspace{1mm}\\
(2)  Every co-semisimple module is a cl.Hilbert module. In fact,
by \cite[Proposition 23.1]{Wis}, \indent  an $R$-module $M$ is
co-semisimple if and only if  every proper submodule of $M$ is
  an  \indent intersection of maximal submodules.\\
(3) Every Artinian $R$-module $M$ is a cl.Hilbert $R$-module (see
Proposition 2.17 (3)).\\
 (4) $\Bbb{Q}$ is not cl.Hilbert $\Bbb{Z}$-module. In general,
let $R$ be an integral domain and $K$ be the \indent quotient
field of $R$. If $K \neq R$ (i.e., $R$ is not a field), then $K$
is not a cl.Hilbert $R$- \indent  module. The zero submodule of
$K$ is a classical prime submodule, but $K$ doesn't \indent have
any maximal $R$-submodule (Let $N$ be a maximal $R$-submodule of
$K$.  Then  \indent ${\rm{Ann}}(K/N)=(0)$ and since $K/N$ is a
simple $R$-module, $(0)$ is  a maximal ideal  of $R$, \indent
i.e., $R$
is field, a contradiction).\\
 (5) Let $R$ be a ring with  dim$(R)=0$. Then every $R$-module is
 cl.Hilbert (see Theorem \indent 2.17 (1)).\\
 (6)  Let $R$  be a Dedekind domain with $J(R)=(0)$. Then every finitely generated  $R$-module \indent is
 cl.Hilbert (see Theorem  3.7).\\
 (7)  Let $R$  be a one-dimensional domain. Then  every torsion $R$-module is   cl.Hilbert (see \indent Theorem  2.17 (2)).

The following two evident lemmas   offer several characterizations
  of classical prime submodules and  prime submodules respectively (see \cite[Propositions 2.1 and 2.2]{BN}
  and also, \cite[Proposition 1.1]{B3}).

\noindent{\bf Lemma 2.3.}  {\it Let $M$ be an $R$-module. For a
submodule $P< M$,
   the following statements
  are equivalent:\vspace{2mm}\\
  $(1)$ $P$ is  classical prime.\vspace{1mm}\\
$(2)$  For every  $0\neq\bar{m}\in M/P$, $(0:R\bar{m})$ is a prime ideal.\vspace{1mm}\\
$(3)$  $\{(0:R\bar{m}) | \ 0\neq\bar{m}\in M/P\}$ is a
chain (linearly ordered set) of prime  ideals.\vspace{1mm}\\
 $(4)$  $(P:M)$ is  a prime ideal, and  $\{(0:R\bar{m})
| \ 0\neq\bar{m}\in M/P\}$ is a chain  of prime  ideals.}

\noindent{\bf Lemma 2.4.}  {\it Let $M$ be an $R$-module. For a
submodule  $P<M$,
   the following statements
  are equivalent:\vspace{2mm}\\
  $(1)$  $P$ is  prime.\vspace{1mm}\\
$(2)$   For every  $0\neq\bar{m}\in M/P$, $(0:R\bar{m})$ is a prime ideal and $(0:R\bar{m})=(P:M)$.\vspace{1mm}\\
$(3)$  $(P:M)$ is  a prime ideal and the set
$\{(0:R\bar{m}):~0\neq\bar{m}\in M/P\}$ is a  singleton.}

In \cite[Theorem 4]{G}, it is shown that a ring $R$ is a Hilbert
ring  if and only if every non-maximal  prime ideal of $R$ is an
intersection of properly larger prime ideals. Next we give a
generalization of this fact to modules.

\noindent{\bf Theorem  2.5.} {\it An $R$-module $M$ is a
cl.Hilbert module if and only if every non-maximal classical prime
submodule of $M$ is an intersection of properly larger classical
prime submodules.}

\noindent{\bf Proof.} If $M$ is a cl.Hilbert module, the given
property certainly holds (since maximal submodules are classical
prime). For the converse, suppose that $N$ is a classical prime
submodule that is not a maximal submodule. Let $m\in M\setminus
N$. Form the set of all classical prime submodule which contain
$N$ but not $m$. This set contains $N$. By Zorn's Lemma, let $K$
be maximal in this set. $K$ must be a maximal submodule.
Otherwise, $K$ is the intersection of properly larger classical
prime submodules. Since $K$ is maximal in the above set of prime
submodules, all properly larger prime submodules must contain $m$.
It would follow from this that $m$ is in $K$. Because this is not
the case, we may conclude that $K$ is indeed a maximal submodule.
We have therefore proved that the intersection of the maximal
submodules which contain $N$ is $N$ itself, and so $M$ is a
cl.Hilbert module.~~$\square$

Let $M$ be an $R$-module and $K\leq M$. One can easily show that a
proper submodule $P$ of $M$ with $K\subseteq P$ is a classical
prime (resp., maximal) submodule of $M$ if and only if $P/K$ is a
classical prime (resp., maximal) submodule of the factor module
$M/K$. The following proposition follows immediately from this
observation.

 \noindent{\bf Proposition 2.6.} {\it Any
homomorphic image of a cl.Hilbert module is a cl.Hilbert module.}

Minimal classical prime submodules are defined in a natural way.
It is clear that whenever $\{P_i\}_{i\in I}$ is a chain of
classical  prime submodules of an $R$-module $M$, then
$\bigcap_{i\in I}P_i$ is always a classical prime submodule. Thus
by Zorn's lemma each classical  prime submodule of $M$ contains a
minimal one (see also \cite[Section 5]{BK}, for more details).

\noindent{\bf Corollary 2.7.} {\it  Let $R$ be a   ring and $M$ be
an  $R$-module. Then the following statements are
equivalent:}\vspace{2mm}\\
(1) {\it  $M$ is a cl.Hilbert $R$-module.}\vspace{1mm}\\
 (2)  {\it  $M/N$ is a cl.Hilbert $R$-module for each
 submodule $N$ of $M$.}\vspace{1mm}\\
(3)  {\it  $M/N$ is a cl.Hilbert $R$-module for each  minimal
classical prime submodule $N$ of  $M$.}

\noindent{\bf Proof.} $(1)\Rightarrow(2)$ is by Proposition 2.6. \\
\indent $(2)\Rightarrow(3)$ is clear.\\
 \indent $(3)\Rightarrow(1).$  Let $P$ be a classical prime
submodule of $M$. Then there is a minimal classical prime $P_0$ of
$M$ contained in $P$. Therefore $P/P_0$ is an intersection of
maximal submodules of $M/P_0$. It follows that $P$ is an
intersection of maximal submodules of $M$.~$\square$

Also by Proposition 2.6, we have the following corollary.

\noindent{\bf Corollary 2.8.}  {\it Let $R$ be a   ring and
$\{M_i\}_{i\in I}$ be a collection of $R$-modules. If
$\bigoplus_{i\in I}M_{i}$ is a cl.Hilbert module, then each
$M_{i}$ $(i\in I)$ is a cl.Hilbert module.}

The next example shows that the converse of Corollary 2.8,
 is not true in general (even if the index set $I$ is finite and each $M_i$ is a finitely generated module).

\noindent{\bf Example 2.9.} Let $R={\Bbb{Z}}[x]$ and $M_1=M_2=R$.
Since $R$ is a Hilbert ring, $M_1$, $M_2$ are cl.Hilbert (Hilbert)
$R$-modules, but by Example 2.1,  $M=M_1\oplus M_2$ is not a
cl.Hilbert $R$-module.

For what follows, We will need the following evident lemma.

  \noindent{\bf Lemma  2.10.} {\it Let $M$ be an $R$-module and  let $I$ be an ideal of $R$
  such that $I \subseteq Ann_{R}(M)$. Then $M$ is a cl.Hilbert $R$-module
if and only if $M$ is a cl.Hilbert $(R/I)$-module.}

Recall that for a ring $R$,  the nilradical of $R$,  denoted by
${\rm Nil}(R)$, is the intersection of all prime ideals of $R$.
Also, for an $R$-module $M$, the radical of $M$, denoted by ${\rm
Rad}_R(M)$, is the intersection of all maximal submodules  of $M$
(if $M$ has no any maximal submodule, then ${\rm Rad}_R(M):=M$).

\noindent{\bf Proposition 2.11.} {\it Let $M$ be an $R$-module.
Then the following statements are equivalent:}\vspace{2mm}\\
\noindent(1) {\it $M$ is cl.Hilbert $R$-module.}\\
\noindent(2) {\it $M/(Nil(R)M)$ is a cl.Hilbert $R$-module}.\\
\noindent(3) {\it $M/(Nil(R)M)$ is a cl.Hilbert
$(R/Nil(R))$-module}.

\noindent{\bf Proof.} (1)$\Rightarrow$(2) is  by Corollary  2.7.\\
\indent (2)$\Rightarrow$(3) is clear by Lemma 2.10.\\
\indent (3)$\Rightarrow$(1). Suppose $P$ is a classical prime
submodule of the $R$-module $M$. Then $(P :M)={\cal{P}}$ is a
prime ideal of $R$ by (4) of lemma 2.3. Thus ${\cal{P}}M\subseteq
P$ and so ${\rm Nil}(R)M\subseteq P$. Now it is clear that $P/{\rm
Nil}(R)M$ is a classical prime submodule of
 $M/{\rm Nil}(R)M$  as an $(R/{\rm Nil}(R))$-module. By our hypothesis we have
  $P/{\rm Nil}(R)M=\bigcap_{i\in I}(M_{i}/{\rm Nil}(R)M)$ where  each $M_{i}/{\rm Nil}(R)M$ is a maximal submodule of $M/{\rm Nil}(R)M$. Hence $P= \bigcap
_{i\in I}M_{i}$, where each  $M_{i}$ is a  maximal submodule of
$M$.~$\square$

\noindent{\bf Proposition 2.12.} {\it Let $M$ be an $R$-module.
Then the following  statements are equivalent:}\vspace{2mm}\\
(1) {\it $M$ is a cl.Hilbert $R$-module.}\vspace{1mm}\\
(2) {\it $M/P$ is a cl.Hilbert $(R/{\cal{P}})$-module  for each
classical prime submodule $P$ of  $M$ with \indent  ${\cal{P}}=(P:M)$.}\vspace{1mm}\\
(3) {\it $Rad_{R}(M/P)=0$  for each classical prime submodule $P$
of $M$}.

\noindent{\bf Proof.} $(1)\Rightarrow (2)$.  Let $P$ be a
classical prime $R$-submodule of $M$ with ${\cal P}=(P:M)$. Then
by Corollary  2.7, $M/P$ is a cl.Hilbert $R$-module. Since ${\cal
P}={\rm{Ann}}(M/P)$, Lemma 2.10 completes the proof.\\
\indent $(2)\Rightarrow (3)$. Let $P$ be a classical prime
submodule of $M$ such that $(P:M)=\cal P$. The zero submodule of
the $(R/{\cal P})$-module $M/P$ is classical prime submodule. By
(2), we have ${\rm Rad}_{R/{\cal P}}(M/P)=0$. On the other hand,
${\rm Rad}_{R/{\cal P}}(M/P)={\rm Rad}_{R}(M/P)=0$.\\
\indent $(3)\Rightarrow (1)$ is clear.~$\square$

\noindent{\bf Proposition 2.13.} {\it  Let $R$ be a domain and $M$
be a cl.Hilbert $R$-module. If $N$ is a any submodule of $M$ such
that $M/N$ is a torsion-free $R$-module, then $N$ is also a
cl.Hilbert  $R$-module.}

\noindent{\bf Proof.} Assume that $R$ is a domain and that $M$ is
a cl.Hilbert $R$-module. Suppose that $N<M$ and that $M/N$ is
torsion-free. Suppose further that $P<N$ is a classical prime
submodule of $N$. We will show that $P$ is the intersection of
maximal submodules of $N$.\vspace{-1.5mm}

We first show that $P$ is a classical prime submodule of $M$.
Toward this end, suppose that $rsm\in P$ for some $m\in M$ and
$r,s\in R$. If $m\in N$, then since $P$ is a classical prime
submodule of $N$, we infer that either $rm\in P$ or $sm\in P$.
Thus assume that $m\notin N$. Recall that $rsm\in P\subseteq N$.
Since $M/N$ is torsion-free and $m\notin N$, it follows that $r=0$
or $s=0$. Thus in this case too, either $rm\in P$ or $sm\in P$.
Thus $P$ is a classical prime submodule of $M$.\vspace{-1.5mm}

Since $P$ is a classical prime submodule of $M$, $P=\cap_{i\in
I}M_{i}$, where each $M_{i}$ is a maximal submodule of $M$. For
each $i$, let $P_{i}:=M_{i}\cap N$. Since $P\subseteq N$, it is
easy to see that $P=\cap_{i\in I}P_{i}$. Further, we may assume
without loss of generality (by discarding all $P_{i}$ containing
$N$, if any) that each $P_{i}$ is properly contained  in $N$. Now
let $i\in I$ be arbitrary. To complete the proof, it suffices to
show that $P_{i}$ is a maximal submodule of $N$. Thus suppose that
$m\in N\setminus P_{i}$. We will show that $(P_{i}, m)=N$. Thus
$m\notin M_{i}$. Since $M_{i}$ is a maximal submodule of $M$, we
have $(M_{i}, m)=M$. Let $x\in N$ be arbitrary (we will show that
$x\in (P_{i}, m)$). Since $M=(M_{i}, m)$, $x=m_{i}+rm$ for some
$m_{i}\in M_{i}$ and $r\in R$. Since $x\in N$ and $m\in N$, we
conclude that $m_{i}\in N$. Thus $m_{i}\in P_{i}$, and it follows
that $x\in (P_{i},m)$. We have shown that $(P_{i}, m)=N$, and this
prove that $P_{i}$ is a maximal submodule of $N$.~$\square$

Recall that a submodule $N$ of an $R$-module $M$ is called {\it
pure} if $IN=N\cap IM$, for every ideal $I$ of $R$. Next, we
easily obtain the following corollary.

\noindent{\bf Corollary 2.14.} {\it Let $R$ be a domain and $M$ be a cl.Hilbert $R$-module. Then the following hold:}\vspace{2mm}\\
(1) {\it  If $T(M)$ is the torsion submodule of $M$, then $T(M)$ is a cl.Hilbert $R$-module.}\vspace{1mm}\\
 (2) {\it If $M$ is torsion-free and $N$ is a pure submodule of $M$, then $N$ is a cl.Hilbert $R$-module.}

\noindent{\bf Proof.} (1) follows immediately from Proposition
2.13. As fore (2), suppose that $N$ is a pure submodule of the
torsion-free cl.Hilbert module $M$. By Proposition 2.13, it
suffices to show that if $m\in M\setminus N$ and $r\in R$ with
$rm\in N$, then $r=0$. So suppose that $m\in M\setminus N$ and
$rm\in N$. Since $N$ is pure, $rM\cap N=rN$. Thus $rm\in rN$, and
there is some $n\in N$ such that $rm=rn$. But then $r(m-n)=0$.
Since $m\notin N$, we see that $m-n\neq 0$. As $M$ is
torsion-free, we conclude that $r=0$. This completes the
proof.~$\square$

We have not found any examples of a cl.Hilbert  module $M$ with a
 submodule $N$ that it is not a  cl.Hilbert  module. Thus an interesting
question is:

\noindent{\bf  Question 2.15.} {\it Is every submodule of a
cl.Hilbert module itself a cl.Hilbert  module?}

Next, we show  that several large classes of modules are classical
Hilbert. We will make use of the following lemma.

\noindent{\bf Lemma 2.16.} {\it Let $R$ be a ring and let $M$ be a
$R$-module. Suppose that $P$ is a classical prime submodule of
$M$. If the set $\{Ann_{R}(m):~0\neq {\bar{m}}\in M/P\}$ consists
only of maximal ideals of $R$, then $P$ is the intersection of
maximal submodules of $M$.}

\noindent{\bf Proof.} Assume that $M$ is an $R$-module and that
$P$ is a classical prime submodule of $M$. Suppose further that
the set  $\{{\rm Ann}_{R}(\bar{m}):0\neq {\bar{m}}\in M/P\}$
consist only of maximal ideals of $R$. By $(4)$ of Lemma 2.3, the
set $\{{\rm Ann}_{R}(\bar{m}):0\neq {\bar{m}}\in M/P\}$ is a
chain. By assumption $\{{\rm Ann}_{R}(\bar{m}):0\neq {\bar{m}}\in
M/P\}$ consist of only maximal ideals of $R$. It follows that
$\{{\rm Ann}_{R}(\bar{m}):0\neq {\bar{m}}\in M/P\}$ is a
singleton, say $\{J\}$. But then ${\rm Ann}_{R}(M/P)=J$, and $M/P$
is naturally a vector space over the field $R/J$. As an
$R/J$-vector space, it is easy to see that the intersection of all
maximal submodules of $M/P$ is also $\{0\}$. Thus $P$ is an
intersection of maximal submodules of $M$.~$\square$

\noindent{\bf Theorem 2.17.} {\it Let $R$ be a ring, and let $M$
be a $R$-module. Then the following hold:}\vspace{2mm}\\
(1) {\it  If $R$ is zero-dimensional, then  $M$ is a cl.Hilbert module.}\vspace{1mm}\\
(2) {\it If $R$ is one-dimensional domain and $M$ is torsion, then
$M$ is a cl.Hilbert module.}\vspace{1mm}\\
(3) {\it If $M$ is Artinian, then $M$ is a cl.Hilbert module.}

\noindent{\bf Proof.} Let $R$ be a ring and let $M$ be an
$R$-module. Suppose that $P$ is a classical prime submodule of $M$
and let $S:=\{{\rm Ann}_{R}(\bar{m}):0\neq {\bar{m}}\in M/P\}$. By
Lemma 2.3, each ${\rm Ann}_{R}(\bar{m})$ is a prime ideal of $R$.
It suffices by Lemma 2.16 to show that if any of the condition in
$(1)-(3)$ hold, then each
${\rm Ann}_{R}(\bar{m})~(0\neq {\bar{m}}\in M/P)$ is a maximal ideal of $R$.\vspace{2mm}\\
\indent (1)  Suppose that $R$ is zero-dimensional. Then as each
${\rm Ann}_{R}(\bar{m})$ is prime, it follows that each
  ${\rm Ann}_{R}(\bar{m})$ is maximal.\vspace{1mm}\\
\indent(2) Assume now that $R$ is a one-dimensional domain and
that $M$ is torsion. Then of course $M/P$ is also torsion. It
follows that
each  ${\rm Ann}_{R}(\bar{m})$ is a nonzero prime ideal of $R$, hence maximal.\vspace{1mm}\\
\indent(3)  Suppose now that $M$ is Artinian, and let $0\neq
{\bar{m}}\in M/P$ arbitrary. Not that $M/P$is Artinian, and hence
also $R\bar{m}$ is Artinian. But $R\bar{m}\cong R/{\rm
Ann}(\bar{m})$, whence $R/{\rm Ann}(\bar{m})$ is an Artinian ring.
Since ${\rm Ann}(m)$ is prime, we see that $R/{\rm Ann}(\bar{m})$
is an Artinian domain, whence a field. Thus ${\rm Ann}(\bar{m})$
is a maximal ideal of $R$.~$\square$

We conclude this section by showing that rings over which all
modules are classical  Hilbert are abundant.

\noindent{\bf Theorem 2.18.} {\it Let $R$ be a  ring. Then the following statements are equivalent:}\vspace{2mm}\\
(1) {\it Every $R$-module is a cl.Hilbert module.}\vspace{1mm}\\
 (2) {\it Every $R$-module is a Hilbert module.}\vspace{1mm}\\
 (3) {\it $dim(R)=0$.}

\noindent{\bf Proof.} $(1)\Rightarrow(2)$ is clear since every prime submodule is classical prime.\\
 \indent $(2)\Rightarrow(3)$. Assume that
every $R$-module is a Hilbert module. Let ${\cal{P}}$ be a prime
ideal of $R$ and let $Q$ be the field of fractions of ${\bar
R}:=R/{\cal{P}}$. Then $(0)< Q$ is a prime ${\bar R}$-submodule.
It follows that $(0)< Q$ is  also a  prime $R$-submodule. If
$Q\neq {\bar {R}}$, then ${\cal{P}}$ is not a maximal  ideal of
$R$ and, one can easily see that $Q$ has no maximal
$R$-submodules, that is a contradiction. Therefore, $Q={\bar
{R}}$, i.e., $P$ is a maximal ideal of $R$ and so ${\rm{dim}}(R)=0$.\\
\indent $(3)\Rightarrow(1)$ is by Theorem 2.17 (1).~$\square$

\section{Rings over which all finitely generated modules are  classical  Hilbert}

 In this section we will  characterize all rings $R$ over which
every  finitely generated  $R$-module is a cl.Hilbert module.

\noindent{\bf Remark  3.1.} Let $R$ be a   ring. Then every
finitely generated $R$-module is a  Hilbert module if and only if
$R$ is a Hilbert ring (see \cite[Proposition 2.9]{MSH}). The
Example 2.1 in Section 2 shows that a finitely generated module
over a Hilbert ring  $R$ need not be a cl.Hilbert $R$-module.  In
fact, in Example 2.1,  it is shown that for the Hilbert ring
${\Bbb{Z}}[x]$ the free ${\Bbb{Z}}[x]$-module
${\Bbb{Z}}[x]\oplus{\Bbb{Z}}[x]$ is not a cl.Hilbert module.

We recall that a {\it Dedekind domain}  is an integral domain $R$
in which every proper  ideal of $R$  is the product of a finite
number of prime ideals. Also, a {\it discrete valuation ring}  is
a principal ideal domain that has exactly one nonzero prime ideal.
A domain  $R$ is a  Dedekind domain if and only if $R$ is
Noetherian and for every nonzero prime ideal ${\cal{P}}$ of $R$,
the localization $R_{\cal{P}}$ of $R$ at ${\cal{P}}$ is a discrete
valuation ring; see for instance, Hungerford \cite[Theorem
6.10]{HU}. Also, it is well-known that a Noetherian local domain
$R$  with maximal ideal ${\cal{M}}$ is a discrete valuation ring
if and only if $R$ is a principal ideal domain, if and only if
${\cal{M}}$  is principal. Thus  we conclude that a Noetherian
domain  $R$ is a Dedekind domain if and only if for every maximal
ideal ${\cal{M}}$ of $R$, the maximal ideal of the localization
$R_{\cal{M}}$ of $R$ at ${\cal{M}}$ is a principal ideal.

We need the following two lemmas.

 \noindent{\bf Lemma 3.2.}
\cite[Lemma 3.3]{B5} {\it Let $R$ be a Dedekind domain. Then every
classical prime submodule of any module is  an intersection of
prime submodules.}

\noindent{\bf Lemma 3.3.} \cite[Proposition 2.4]{B5} {\it Suppose
that $M$ is a Noeitherian module over a  ring $R$. Then the
   following statements are equivalent:}\vspace{1mm}\\
   (1) {\it  Every  classical prime submodule of $M$ is  an intersection of prime
   submodules.}\vspace{1mm}\\
   (2) {\it  For every maximal ideal ${\mathcal{M}}$ of $R$,  every classical  prime submodule of $M_{\mathcal{M}}$ as   an
   $R_{\mathcal{M}}$- \indent module  is  an intersection of prime
   submodules.}

In \cite[Theorem 3.5]{B5}, it is shown that if $R$ is  a
commutative Noetherian  domain, then  every  classical prime
submodule of $M$ is  an intersection of prime
   submodules if and only if   $R$
   is a Dedekind domain.  In what follows, we show that if even  every
classical prime submodule of the free module $R\oplus R$ is  an
intersection of prime submodules, then  $R$
   is a Dedekind domain.

   \noindent{\bf Theorem 3.4.} {\it Let $R$ be a   Noetherian  domain. Then the following statements are equivalent.}\vspace{2mm}\\
    (1) {\it Every  classical prime submodule of any module is  an intersection of prime
     submodules.}\vspace{1mm}\\
   (2) {\it Every  classical prime submodule of each  finitely generated  module is  an intersection \indent of prime
   submodules.}\vspace{1mm}\\
   (3) {\it Every  classical prime submodule of the free module $R\oplus R$ is  an  intersection    of prime
  \indent  submodules.}\vspace{1mm}\\
    (4) {\it  $R$    is a Dedekind domain.}

   \noindent{\bf Proof.} $(1)\Rightarrow(2)\Rightarrow(3)$ is clear.\\
   \indent $(3)\Rightarrow(4)$. We can assume that  $R$ is not a field. Then ${\rm{dim}}(R)\geq
   1$. Since $R$ is a Noetherian domain, it suffices to show that for every maximal ideal
${\cal{M}}$ of $R$, the maximal ideal  ${\cal{M}}^e$ of the
localization $R_{\cal{M}}$ of $R$ at ${\cal{M}}$ is a principal
ideal. Let $\mathcal{M}$ be the maximal ideal of $R$. By $(3)$ and
Lemma 3.3, every classical $R_{\cal{M}}$-submodule of the free
$R_{\cal{M}}$-module $R_{\cal{M}}\oplus R_{\cal{M}}$ is  an
intersection of prime submodules. Thus we may assume that $R$ is a
local domain. Choose $a\in {\cal{M}}\setminus {\cal{M}}^2$. If
   ${\mathcal{M}}=Ra$, then we are done. Suppose not. Then we
   can choose $b\in {\mathcal{M}}\setminus Ra$. As a $a\in {\mathcal{M}}\setminus
   {\mathcal{M}}^2$, $a\in{\mathcal{M}}\setminus Rb$. It follows that
   $\Lambda(a,b)=\{(x,y)\in R\oplus R  :  xb=ya \}\subseteq
   \mathcal{M}\oplus\mathcal{M}$. It is easily checked that
   $\Lambda(a,b)$ is a prime submodule of $R\oplus R$. Now we
   claim that ${\mathcal{M}}\Lambda(a,b)$ is a classical prime
   submodule of $R\oplus R$. To see this, let $rs(x,y)\in {\mathcal{M}}\Lambda(a,b)$,
    where $(x,y)\in R\oplus R\setminus {\mathcal{M}}\Lambda(a,b)$ and $r, s \in
    R\setminus(0)$. Therefore  $rs(x,y)\in\Lambda(a,b)$, which implies that   either we have  $(x,y)\in\Lambda(a,b)$ or
    $rs(R\oplus R)\subseteq \Lambda(a,b)$.  But if  $rs(R\oplus R)\subseteq \Lambda(a,b)$, then
    $rs(1,1)\in\Lambda(a,b)$ and we must have $a=b$, which is a  contradiction.
    Thus we must have, $(x,y)\in \Lambda(a,b)$ and therefore
    $s(x,y)\in{\mathcal{M}}\Lambda(a,b)$, which means that ${\mathcal{M}}\Lambda(a,b)$  is a classical  prime submodule of
    $R\oplus R$. Now by our hypothesis ${\mathcal{M}}\Lambda(a,b)$ is  an
intersection of prime submodules of $R\oplus R$. Let $P$ be a
prime submodule of $R\oplus R$ that contains
${\mathcal{M}}\Lambda(a,b)$. We have
     ${\mathcal{M}}(R\oplus R)\subseteq P$ or  $\Lambda(a,b)\subseteq P$. In any case, $\Lambda(a,b)\subseteq P$ (since $\Lambda(a,b)\subseteq
   {\cal{M}}\oplus{\cal{M}}={\cal{M}}(R\oplus R)$).  It follows that  ${\mathcal{M}}\Lambda(a,b)=\Lambda(a,b)$.
      By Nakayama's Lemma,  $\Lambda(a,b)=(0)$ which contradicts $(a,b)\in\Lambda(a,b)$. Therefore, ${\mathcal{M}}=Ra$
       and so  $R$ is a Dedekind domain.\\
      \indent  $(4)\Rightarrow (1)$ is by Lemma 3.2.~$\square$

We also need the following  lemma.

  \noindent{\bf Lemma 3.5.} {\it Let $R$ be a Dedekind domain with $J(R)=(0)$.
   Then every finitely generated $R$-module is a
cl.Hilbert module.}

\noindent{\bf Proof.} Let $R$ be a Dedekind domain with $J(R)=(0)$
and let $M$ be a finitely generated $R$-module. Clearly $R$ is a
Hilbert ring and so by \cite[Proposition 2.9]{MSH}, $M$ is a
Hilbert module. Since $R$ is a Dedekind domain,  by Lemma 3.2,
every classical prime submodule of $M$ is  an intersection of
prime submodules of $M$. Thus every classical prime submodule of
$M$ is an intersection of maximal submodules of $M$, i.e., $M$ is
a cl.Hilbert module. $\square$

Now we are in ready to characterize those commutative Noetherian
domains $R$ over which all finitely generated $R$-modules are
cl.Hilbert.

\noindent{\bf Theorem  3.6.} {\it Let $R$ be a Noetherian  domain.
Then the following statements are
equivalent:}\vspace{2mm}\\
(1) {\it Every finitely generated $R$-module is a
cl.Hilbert module.}\vspace{1mm}\\
(2) {\it The free $R$-module $R\oplus R$ is a
cl.Hilbert module.}\vspace{1mm}\\
(3) {\it $R$ is both  a Hilbert ring and a  Dedekind domain.}\vspace{1mm}\\
(4) {\it $R$ is a  Dedekind domain with $J(R)=0$.}\vspace{1mm}\\
(5) {\it $R$ is either a field or a  Dedekind domain with infinity
many  maximal ideals.}

\noindent{\bf Proof.} $(1)\Rightarrow(2)$ is clear.\\
\indent $(2)\Rightarrow(3)$. Since the free $R$-module $R\oplus R$
is a cl.Hilbert module, it follows   by Corollary 2.7 that the
$R$-module $R$ also a cl.Hilbert module, i.e., $R$ is a Hilbert
ring. Since every classical prime submodule of the free $R$-module
$R\oplus R$ is an intersection of maximal (prime)  submodules, it
follows by Theorem 3.4, $R$
is a Dedekind domain.\\
\indent $(3)\Rightarrow(4)$. Since $R$ is a Hilbert domain, $(0)$
is an
intersection of maximal ideals, i.e., $J(R)=(0)$.\\
\indent $(4)\Rightarrow(1)$ is by Lemma 3.5. \\
\indent $(4)\Rightarrow(5)$. Suppose that $R$ is not a field.
Since $R$ is a domain with $J(R)=(0)$, we conclude that the set of
maximal ideals of $R$
is infinite.\\
\indent $(5)\Rightarrow(4)$. Suppose to contrary that  $J(R)\neq
(0)$. Then ${\rm{dim}}(R/J(R))=0$  and since $R$ is  Noetherian,
we conclude that $R/J(R)$ is an Artinian ring with infinity many
maximal ideals, a contradiction.~$\square$

Finally,  we  characterize    Noetherian rings $R$ over which all
finitely generated $R$-modules are cl.Hilbert.

\noindent{\bf Theorem  3.7.} {\it Let $R$ be a
 ring. Consider the following statements.}\vspace{2mm}\\
(1) {\it Every finitely generated $R$-module is a
cl.Hilbert module.}\vspace{1mm}\\
(2) {\it Every finitely generated $R/{\cal{P}}$-module is a
cl.Hilbert module  for each minimal   prime \indent ideal ${\cal{P}}$ of $R$.}\vspace{1mm}\\
(3) {\it The free $R$-module $R\oplus R$ is a
cl.Hilbert module.}\vspace{1mm}\\
(4) {\it The  free $R/{\cal{P}}$-module  $R/{\cal{P}}\oplus R/{\cal{P}}$ is  a
cl.Hilbert module for each minimal    prime ideal \indent ${\cal{P}}$ of $R$.}\vspace{1mm}\\
(5) {\it $R$ is a Hilbert ring and for each minimal prime ideal
${\cal{P}}$ of $R$, the ring $R/{\cal{P}}$ is    a \indent
Dedekind
 domain.}\vspace{1mm}\\
{\it Then $(1)\Leftrightarrow(2)\Rightarrow(3)\Leftrightarrow(4)$
and $(5)\Rightarrow(1)$. When $R$ is a Noetherian ring, all the
five statements are equivalent.}

\noindent{\bf Proof.} $(1)\Rightarrow(2)$. Since every $R/{\cal{P}}$-module is
an $R$-module by $rm:=(r+{\cal{P}})m$, the proof is clear.  \\
\indent $(2)\Rightarrow(1)$. Let  $M$ be a finitely generated
$R$-module and $P$ be  a classical prime submodule of $M$. Then
${\cal{P}}=(P:M)$ is a prime ideal of $R$. Suppose that
${\cal{P}}_0\subseteq {\cal{P}}$ is a minimal prime ideal of $R$.
Then $M/N$ is a classical  $R/{\cal{P}}_0$-module and so  by our
hypothesis $M/N$ is a cl.Hilbert $R/{\cal{P}}_0$-module. Thus the
zero submodule of $M/N$ is an intersection of maximal
$R/{\cal{P}}_0$-submodules of $M/N$. It follows that $N$ is  an
intersection of maximal $R$-submodules of $M$. Thus $M$ is a cl.Hilbert $R$-module.\\
\indent $(2)\Rightarrow(3)$ is clear.\\
\indent $(3)\Leftrightarrow(4)$ is similar to the proof of
$(1)\Leftrightarrow(2)$.\\
\indent $(5)\Rightarrow(1)$. Let  $M$ be a finitely generated
$R$-module and $P$ be  a classical prime submodule of $M$. Then
${\cal{P}}=(P:M)$ is a prime ideal of $R$. Suppose that
${\cal{P}}_0\subseteq {\cal{P}}$ is a minimal prime ideal of $R$.
Then $M/N$ is a classical  $R/{\cal{P}}_0$-module. Since $R$ is a
Hilbert ring,  $R/{\cal{P}}_0$ is also a Hilbert ring and by our
hypothesis  $R/{\cal{P_0}}$ is a  Dedekind domain. Thus by Lemma
3.5, $M/P$ is a cl.Hilbert $R/{\cal{P}}_0$-module. Thus the zero
submodule of $M/N$ is an intersection of maximal
$R/{\cal{P}}_0$-submodules of $M/N$. It follows that $N$ is  an
intersection of maximal $R$-submodules of $M$. Thus $M$ is a
cl.Hilbert $R$-module.\\
\indent For the proof of the second statement,  we show that
$(3)\Rightarrow(5)$. Assuming that $R$  is a Noetherian ring.
Since the  free $R$-module  $R\oplus R$ is a cl.Hilbert module, we
conclude that $R$ is a Hilbert ring. It follows that for each
minimal  prime ideal ${\cal{P}}$ of $R$ the ring $R/{\cal{P}}$ is
a Hilbert ring and also a Noetherian domain.  On the other hand,
by $(3)\Leftrightarrow (4)$, the  free $R/{\cal{P}}$-module
  $R/{\cal{P}}\oplus R/{\cal{P}}$ is  a cl.Hilbert module. Thus by Theorem 3.4,
   the ring $R/{\cal{P}}$ is a Dedekind domain.~$\square$

\end{document}